# Modularity of fibres in rigid local systems

By Henri Darmon

## 1. Introduction

Let $K$ be a totally real field embedded in a fixed algebraic closure $\overline{K}$, and write $G_K := \mathrm{Gal}(\overline{K}/K)$ for its absolute Galois group. Fix a prime $\ell \neq 2$, and consider an odd two-dimensional Galois representation

$$\rho : G_K \longrightarrow \mathbf{GL}_2(E),$$

where $E$ is either a finite field of characteristic $\ell$ or a finite extension of $\mathbb{Q}_\ell$. Assume that the restrictions of $\rho$ to the inertia groups at the primes of $K$ above $\ell$ are *potentially semistable* in the sense of [FM].

The representation $\rho$ is called *modular* if it is associated to a Hilbert modular form on $\mathbf{GL}_2(K)$, as is explained, for example, in [W1] and [W2]. Fontaine and Mazur [FM] conjectured that this is always the case. Significant progress on this conjecture was achieved [W3] by proving particular instances of the following "lifting conjecture":

CONJECTURE 1.1. *Suppose that $\ell$ is odd and that the residual representation $\bar\rho$ attached to $\rho$ is modular. Then $\rho$ itself is modular.*

Conjecture 1.1 is proved in [W3] and [TW] when $K = \mathbb{Q}$ and the restriction of $\rho$ to the decomposition groups at the primes above $\ell$ are *semistable* in the sense of [DDT, §2.4]. This is enough (using the primes $\ell = 3$ and $5$) to establish the Shimura-Taniyama conjecture for semistable elliptic curves, thanks to a key result of Langlands and Tunnell. Progressively stronger cases of Conjecture 1.1 were subsequently proved by [Di], [CDT], [Fu], and [SW1]; in [SW1], Skinner and Wiles obtain quite general results in the context where $K$ is any totally real field, the principal assumption being that $\rho$ is *ordinary* at the primes above $\ell$.

In this note we consider Galois representations which occur in "rigid families", and establish their modularity under Conjecture 1.1. This implies the modularity (over suitable real abelian extensions) of the Galois representations occurring in the cohomology of the curves

$$y^n = x^a(x-1)^b(x-t)^c, \qquad t \in \mathbb{Q},$$

whose periods as a function of the parameter $t$ are values of classical hypergeometric functions.

To state the main result precisely, denote by $K(t)$ the field of rational functions in the indeterminate $t$, and let

$$\varrho : G_{K(t)} \longrightarrow \mathbf{GL}_2(E)$$

be a two-dimensional Galois representation. For $x \in \mathbb{P}^1(\bar{K})$, viewed as a place of $K(t)$, let $D_x \subset G_{K(t)}$ be a decomposition group at $x$, and write $I_x = \hat{\mathbb{Z}}(1)$ for its inertia subgroup. One says that $\varrho$ is *unramified* at $x$ if its restriction to $I_x$ is trivial. If, in addition, $x$ belongs to $\mathbb{P}^1(K)$, the restriction of $\varrho$ to $D_x$ factors through $D_x/I_x = G_K$, giving rise to a Galois representation

$$\varrho[x] : G_K \longrightarrow \mathbf{GL}_2(E),$$

which can be thought of as the specialization of $\varrho$ at $t = x$.

Let $\varrho^{\mathrm{geom}}$ be the restriction of $\varrho$ to the subgroup

$$G^{\mathrm{geom}} := \mathrm{Gal}(\overline{K(t)}/\overline{K}(t)) \subset G_{K(t)}.$$

The representation $\varrho$ is said to be *rigid* if $\varrho^{\mathrm{geom}}$ is unramified outside 0, 1, and $\infty$. (The reason for this terminology will be made clear in the next section; cf. Prop. 2.4.) Choose a topological generator of $\hat{\mathbb{Z}}(1)$ corresponding to a compatible system $(\zeta_n)$ of primitive $n^{\mathrm{th}}$ roots of unity. For $j = 0, 1, \infty$, let $\gamma_j$ be the corresponding generator of $I_j$, and let $\sigma_j = \varrho(\gamma_j) \in \mathbf{GL}_2(E)$. The *monodromy matrices* $\sigma_j$ depend on the choice of decomposition groups $D_j$ but their conjugacy classes in $\mathbf{GL}_2(E)$ are well-defined. We will show (Lemma 2.2) that the semisimplification of $\sigma_j$ has finite order $n_j$. One can then prove (Prop. 2.4) that the "field of definition" $K$ of $\varrho$ necessarily contains the real subfield $K_{n_j} := \mathbb{Q}(\zeta_{n_j})^+$ of the cyclotomic field of $n_j$-th roots of unity. Conversely, $\varrho$ has a twist which extends to a representation of $G_{K_n(t)}$, where $n = n(\varrho)$ is the least common multiple of the $n_j$. Replace $\varrho$ by such a twist, and $K$ by $K_n$. Our main result is then:

THEOREM 1.2. *Let $\varrho$ be a rigid representation, and assume that one of the $\sigma_j$ is unipotent, and that 8 does not divide $n = n(\varrho)$.*

*If Conjecture 1.1 is true, then $\varrho[x]$ arises from a Hilbert modular form over $K_n$, for all $x \in \mathbb{P}^1(\mathbb{Q}) - \{0, 1, \infty\}$.*

*Acknowledgements.* The author thanks Nick Katz for helpful conversations and the ETH in Zürich for its hospitality while this article was written. This research was funded by grants from NSERC and by an Alfred P. Sloan research fellowship.

## 2. Rigid representations

Fix a rigid representation $\varrho$, and keep the notations of the introduction. While the monodromy matrices $\sigma_j$ are only defined up to conjugation, the decomposition groups $D_j$ can be chosen so that the relation

(1) $$\sigma_0 \sigma_1 \sigma_\infty = 1$$

is satisfied (cf. for example [Se1, Th. 6.3.2]). Fix such a choice from now on.

A $2 \times 2$ matrix is called a *reflection* if its eigenvalues are $1$ and $-1$.

LEMMA 2.1. *The matrix $\sigma_j$ is either a reflection or an element of $\mathbf{SL}_2(E)$.*

*Proof.* The conjugacy classes of $\sigma_j$ are *rational* over the real field $K$ in the sense of [Se1, Sec. 7.1]. In particular, $\sigma_j$ is conjugate to $\sigma_j^{-1}$; the result follows.

If one of the $\sigma_j$ is a reflection, then exactly two are, because of the relation $\det(\sigma_0 \sigma_1 \sigma_\infty) = 1$ which follows from (1). In that case, the image of $\varrho^{\text{geom}}$ is a dihedral group. We exclude this case from consideration from now on, and assume that each $\sigma_j$ belongs to $\mathbf{SL}_2(E)$. The matrix $\sigma_j$ is said to be *quasi-unipotent* if its minimal polynomial has a double root.

LEMMA 2.2. *The $\sigma_j$ are either quasi-unipotent or of finite order.*

*Proof.* Let $K^{\text{cyc}} := K(\zeta_\infty)$ be the maximal cyclotomic extension of $K$, and let $\Omega$ be its Galois group, identified with a subgroup of $\hat{\mathbb{Z}}^\times$. Since the conjugacy class of $\sigma_j$ is rational over $K$, the matrix $\sigma_j$ is conjugate to $\sigma_j^\alpha$ for all $\alpha \in \Omega$. But $\Omega$ has finite index in $\hat{\mathbb{Z}}^\times$, and hence the eigenvalues of $\sigma_j$ are roots of unity.

*Definition* 2.3. An *admissible triple* in $\mathbf{SL}_2(E)$) is a triple $(\sigma_0, \sigma_1, \sigma_\infty)$ of elements in $\mathbf{SL}_2(E)$, taken modulo conjugation in $\mathbf{GL}_2(E)$, and satisfying

(a) The semisimplification of $\sigma_j$ has finite order $n_j$;

(b) The group generated by $\sigma_0$, $\sigma_1$, and $\sigma_\infty$ is an irreducible subgroup of $\mathbf{SL}_2(E)$.

(c) $\sigma_0 \sigma_1 \sigma_\infty = 1$.

Let $n = n(\varrho) = \text{lcm}(n_0, n_1, n_\infty)$, as before. The following "rigidity" property justifies the terminology of the introduction.

PROPOSITION 2.4. *Let $(\sigma_0, \sigma_1, \sigma_\infty)$ be an admissible triple in $\mathbf{SL}_2(E)$ with $\sigma_1$ unipotent. Then there exists a rigid representation*

$$\varrho : G_{K_n(t)} \longrightarrow \mathbf{GL}_2(E)$$

*whose monodromy matrix at* $t = j$ *is equal to* $\sigma_j$. *Furthermore, if* $\varrho'$ *is any irreducible rigid representation whose monodromy matrices are conjugate to those of* $\varrho$, *then* $\varrho'$ *is conjugate to* $\varrho \otimes \chi$, *where* $\chi : G_{K_n} \longrightarrow E^\times$ *is a constant central character.*

This follows from Theorems 1 and 2 of [Be]. (See also the discussion in Section 1 of [Da2].)

### 3. Hypergeometric abelian varieties

For the following definition, let $K$ be any real abelian field, and $\mathcal{O}_K$ its ring of integers. (We will also write $\mathcal{O}_n := \mathbb{Z}[\zeta_n + \zeta_n^{-1}]$ to denote the ring of integers of $K_n$.)

*Definition* 3.1. A *hypergeometric abelian variety with multiplications by* $K$ is an abelian scheme $A$ over $(\mathbb{P}^1 - \{0, 1, \infty\})_{/\mathbb{Q}}$ of dimension $[K : \mathbb{Q}]$ equipped with an inclusion
$$\iota : \mathcal{O}_K \hookrightarrow \operatorname{End}_{K(t)}(A)$$
which is compatible with the natural action of $\operatorname{Gal}(K/\mathbb{Q})$ on both sides, and whose associated monodromy representation is irreducible.

Define an *admissible triple* $(\sigma_0, \sigma_1, \sigma_\infty)$ in $\mathbf{SL}_2(\mathcal{O}_K)$ in the obvious way (replacing $E$ by $\mathcal{O}_K$ in Definition 2.3). Given a hypergeometric abelian variety $A$ with multiplications by $K$, one can associate to it an admissible triple $(\sigma_0, \sigma_1, \sigma_\infty)$ in $\mathbf{SL}_2(\mathcal{O}_K)$ by letting $\sigma_j$ be the image of $\gamma_j$ acting on the DeRham cohomology $H^1_{\mathrm{Dr}}(A)$ (viewed as a two-dimensional $K$-vector space). Conversely, given an admissible triple $(\sigma_0, \sigma_1, \sigma_\infty)$ in $\mathbf{SL}_2(\mathcal{O}_K)$, let $n_j$ be the order of the semisimplification of $\sigma_j$ and set $n = \operatorname{lcm}(n_0, n_1, n_\infty)$. One sees that $K$ must contain the fields $K_{n_j}$ generated by the traces of the $\sigma_j$. Assume that $K = K_n$.

PROPOSITION 3.2. *Assume that* $\sigma_1$ *is unipotent. There exists a hypergeometric abelian variety* $A$ *with multiplications by* $K_n$ *whose associated monodromies are* $(\sigma_0, \sigma_1, \sigma_\infty)$. *The isogeny class of this abelian variety depends only on the triple* $(\sigma_0, \sigma_1, \sigma_\infty)$ *(modulo conjugation by* $\mathbf{GL}_2(E)$).

*Proof* (See [Ka, §5.4], or [CW, §3.3]). The hypergeometric abelian varieties are constructed as appropriate quotients of the Jacobians of the curves

$$y^n = x^a(x-1)^b(x-t)^c.$$

*From hypergeometric abelian varieties to rigid representations.* If $A$ is a hypergeometric abelian variety with multiplications by $K$, the $\ell$-adic Tate module
$$T_\ell(A) := \varprojlim A[\ell^k]$$
is a free module of rank two over $\mathcal{O}_K \otimes \mathbb{Z}_\ell$. The natural action of $G_{\mathbb{Q}(t)}$ on this Tate module is semilinear, in the sense that
$$\alpha(s \cdot v) = s^\alpha \cdot \alpha(v), \quad \text{for} \quad \alpha \in G_{\mathbb{Q}(t)}, \quad s \in \mathcal{O}_K \otimes \mathbb{Z}_\ell, \quad v \in T_\ell(A).$$
In particular, if $\varphi$ is a homomorphism from $\mathcal{O}_K$ to $E$, then $T_\ell(A) \otimes_\varphi E$ is a two-dimensional $E$-vector space on which $G_{K(t)}$ acts linearly. It gives rise to a rigid two-dimensional Galois representation $\varrho$ of $G_{K(t)}$, and thus to a family of representations $\varrho[x]$ of $G_K$ for all $x \in K - \{0, 1\}$.

*Definition* 3.3. The hypergeometric abelian variety $A$ is said to be *modular* at $x$ if $\varrho[x]$ is associated to a Hilbert modular form over $K$ with coefficients in $E$, for all choices of $(\varphi, E)$. We say that $A$ is *modular* if it is modular at $x$, for all $x \in \mathbb{Q} - \{0, 1\}$.

*Remark.* The representations $\varrho[x]$ attached to $A$, as $\ell$, $E$, and $\varphi$ vary, form a *compatible system* of $\ell$-adic representations of $G_K$, and hence prove that $A$ is modular at $x$, it suffices to prove that $\varrho[x]$ is modular for a single $E \subset \bar{\mathbb{Q}}_\ell$.

*Examples.*

1. If $\sigma_0, \sigma_1$, and $\sigma_\infty \in \mathbf{SL}_2(\mathbb{Z})$ are quasi-unipotent with eigenvalues 1, 1, and $-1$, then $A$ is isogenous to the Legendre family of elliptic curves
$$y^2 = x(x-1)(x-t).$$
The modularity of $A$ is thus a special case of the Shimura-Taniyama conjecture which was completely established by Wiles [W3].

2. If $\sigma_0$ and $\sigma_\infty$ are of order 4 and 3, respectively, and $\sigma_1$ is unipotent, then $A_{/\mathbb{Q}(t)}$ is isogenous to (a twist of) the universal family of elliptic curves of invariant $j = 1728/(t-1)$. The modularity of $A$ in this case is merely a reformulation of the Shimura-Taniyama conjecture.

3. If $\sigma_0$ and $\sigma_1$ are unipotent and $\sigma_\infty$ is of order $r$ with $r$ an odd prime, then the corresponding hypergeometric abelian variety is the Jacobian of the hyperelliptic curve with real multiplications by $\mathbb{Q}(\zeta_r)^+$ given by the equation
$$y^2 = (x+2)(f(x) + 2 - 4t),$$
where $f(x) = xg(x^2 - 2)$ and $g(x)$ is the characteristic polynomial of $-(\zeta_r + \zeta_r^{-1})$. This curve had already been considered in [TTV], and used in [Da2] to study the generalized Fermat equation $x^p + y^p = z^r$. In the language of [Da2], the mod $p$ representations attached to $A$ are the "even Frey representations" associated to the generalized Fermat equation $x^p + y^p = z^r$.

4. If $\sigma_0$ and $\sigma_1$ are unipotent and $\sigma_\infty$ has order $2r$ with $r$ an odd prime, then $A$ is the Jacobian of the hyperelliptic curve (also used in the study of $x^p + y^p = z^r$)
$$y^2 = f(x) + 2 - 4t.$$

*From rigid representations to hypergeometric abelian varieties.* Let $\varrho$ be a rigid representation of $G_{K_n}(t)$ with unipotent monodromy at $t = 1$, associated to an admissible triple $(\sigma_0, \sigma_1, \sigma_\infty)$ in $\mathbf{SL}_2(E)$. This triple can be lifted to an admissible triple $(\tilde\sigma_0, \tilde\sigma_1, \tilde\sigma_\infty)$ in $\mathbf{SL}_2(\mathcal{O}_n)$, i.e., there is a homomorphism $\varphi : \mathcal{O}_n \longrightarrow E$ such that $\varphi(\tilde\sigma_j) = \sigma_j$, and $\tilde\sigma_1$ is unipotent. Let $A$ be the hypergeometric abelian variety with multiplications by $K_n$ associated to $(\tilde\sigma_0, \tilde\sigma_1, \tilde\sigma_\infty)$ by Proposition 3.2. Then we have:

PROPOSITION 3.4. *The representation $\varrho$ is equivalent to (a twist of) the Galois representation obtained from the action of $G_{K(t)}$ on $T_\ell(A) \otimes_\varphi E$.*

*Proof.* This is a direct consequence of the uniqueness statement of Proposition 2.4, since the representation associated to $T_\ell(A) \otimes_\varphi E$ is a rigid representation associated to the triple $(\sigma_0, \sigma_1, \sigma_\infty)$. □

Thanks to Proposition 3.4, it is enough to show that all hypergeometric abelian varieties with unipotent monodromy at $t = 1$ and $8 \nmid n$ are modular in order to prove Theorem 1.2.

## 4. Congruences

Let $A$ be a hypergeometric abelian variety with multiplication by $K = K_n$, and let $(\sigma_0, \sigma_1, \sigma_\infty)$ be the associated admissible triple in $\mathbf{SL}_2(\mathcal{O}_K)$. Assume that $\sigma_1$ is unipotent, and let $\ell$ be an odd prime which divides $n = \mathrm{lcm}(n_0, n_1, n_\infty)$. For $j = 0, 1, \infty$, let $n_j'$ be the prime-to-$\ell$ part of $n_j$, let $n'$ be the prime-to-$\ell$ part of $n$, and let $K' = \mathbb{Q}(\zeta_{n'})^+$. Choose a prime $\lambda$ of $K$ above $\ell$, and let $\lambda'$ be the unique prime of $K'$ below it. The prime $\lambda'$ is totally ramified in $K/K'$, so that the residue fields of $K$ and $K'$ at $\lambda$ and $\lambda'$ respectively are canonically isomorphic. Let $\mathbb{F}$ be this common residue field. It is equipped with maps $\varphi : \mathcal{O}_K \longrightarrow \mathbb{F}$ and $\varphi' : \mathcal{O}_{K'} \longrightarrow \mathbb{F}$. Let $(\sigma_0', \sigma_1', \sigma_\infty')$ be a lift of $(\varphi(\sigma_0), \varphi(\sigma_1), \varphi(\sigma_\infty))$ to an admissible triple in $\mathbf{SL}_2(\mathcal{O}_{K'})$, and let $A'$ be the abelian variety associated to it by Proposition 3.2.

Because $G_{\mathbb{Q}(t)}$ acts semi-linearly on $A[\ell] \otimes_\varphi \mathbb{F}$ and because $\lambda$ is totally ramified in $K/K'$, the action of $G_{K(t)}$ on this $\mathbb{F}$-vector space extends to a linear action of $G_{K'(t)}$.

THEOREM 4.1. *The $G_{K'(t)}$ representation $A[\ell] \otimes_\varphi \mathbb{F}$ is isomorphic to (a twist of) the representation $A'[\ell] \otimes_{\varphi'} \mathbb{F}$.*

*Proof.* A direct consequence of Proposition 3.4.

## 5. Proof of the main result

THEOREM 5.1. *Let $A$ be a hypergeometric abelian variety with multiplications by $K_n$, and let $(\sigma_0, \sigma_1, \sigma_\infty)$ be the associated admissible triple. Assume that $\sigma_1$ is unipotent, and that 8 does not divide $n$. If Conjecture 1.1 is true, then $A$ is modular.*

*Proof.* The proof is by induction on $d = [K_n : \mathbb{Q}]$. If $d = 1$, then $A$ is an elliptic curve over $\mathbb{Q}(t)$ and the modularity of $A$ follows from the Shimura-Taniyama conjecture, which itself follows from Conjecture 1.1. If $d > 1$, then $n$ is divisible by an odd prime $\ell$, by the assumption that 8 does not divide $n$. Adopting the notation of Section 4, we begin by showing (for a fixed $t = x \in \mathbb{Q}$) that $A[\ell] \otimes_\varphi \mathbb{F}$ is associated to a Hilbert modular form $f_\ell$ over $K$. If $A[\ell] \otimes_\varphi \mathbb{F}$ is a reducible representation of $G_K$, then one may express $f_\ell$ in terms of Eisenstein series. Assume that $A[\ell] \otimes_\varphi \mathbb{F}$ is irreducible. Since $n' < n$ and $d' = [K' : \mathbb{Q}] < d$, the induction hypothesis implies that $A'$ is modular. Hence so is the rigid representation $A'[\ell] \otimes_{\varphi'} \mathbb{F}$; let $f'_\ell$ be the associated Hilbert modular form mod $\ell$ on $\mathbf{GL}_2(K')$. By Theorem 4.1, the $G_{K'}$ module $A[\ell] \otimes_\varphi \mathbb{F}$ is isomorphic to $A'[\ell] \otimes_{\varphi'} \mathbb{F}$, and so corresponds to the same $f'_\ell$. Letting $f_\ell$ be the cyclic base change lift (from $K'$ to $K$) of $f'_\ell$, it follows that the representation $A[\ell] \otimes_\varphi \mathbb{F}$ is modular over $K$. The $\lambda$-adic Tate module $T_\ell(A) \otimes K_\lambda$ is a potentially semistable Galois representation, since it arises from the torsion points of an abelian variety. Hence it is modular, by Conjecture 1.1.

*Remark.* The proof that $A$ is modular involves repeated applications of the lifting Conjecture 1.1, once with each odd prime $\ell$ dividing $n$. In light of the results in [SW1], it might be feasible to prove unconditionally the modularity of $A$ at $t = x$, when $x$ is such that $A$ is *ordinary* at all these primes. There are infinitely many values of $x$ with this property: for example, all the $x$ for which $n$ divides the numerator of $x - 1$.

MCGILL UNIVERSITY, MONTREAL, PQ, CANADA
*E-mail address*: darmon@math.mcgill.edu